\documentclass{article}
\usepackage{amsmath,amssymb,graphicx,subfig,multirow,comment,fancyhdr}
\newtheorem{theorem}{Theorem}[section]

\newcommand{\norm}[1]{\left\|#1\right\|} 

\title{An adaptive iterative/subdivision hybrid algorithm for curve/curve intersection}

\author{Gun Srijuntongsiri\thanks{Sirindhorn International Institute of Technology, Thammasat University, 131 Moo 5, Tiwanont Road, Bangkadi,
Muang, Pathum Thani, 12000, Thailand. Email: gun@siit.tu.ac.th.}}

\begin{document}

\maketitle

\begin{abstract}
The behavior of the iterative/subdivision hybrid algorithm for curve/curve intersection proposed in \cite{srijuntongsiri_cci} depends on the choice of the domain for their convergence test.  Using either too small or too large test domain may cause the test to fail to detect cases where Newton's method in fact converges to a solution, which results in unnecessary additional subdivisions and consequently more computation time.  
We propose a modification to the algorithm to adaptively adjust the test domain size according to what happens during the test of the parent region.  This is in contrast to the original algorithm whose test domain is always a fixed multiple of the input domain under consideration.  Computational results show that the proposed algorithm is slightly more efficient than the original algorithm.


\end{abstract}

\section{Introduction}

We consider the problem of finding all of the intersection points between two space curves, which is typically known as curve/curve intersection (CCI).  This problem is fundamental in geometric modeling and design, computational geometry, and manufacturing applications \cite{hoschek,mortenson}, and has applications in robotics, collision avoidance, manufacturing simulation, and scientific visualization, for example.  CCI is also a basis for solving more general intersection problems, such as finding intersections between two surfaces.  

There are many different approaches to solve CCI such as subdivision methods \cite{whitted,rubi,houghton,patrikalakis1,sederberg_clip}, iterative methods \cite{heath}, genetic algorithms \cite{hawat}, and algorithms based on orthogonal projections \cite{limaiem}.  Srijuntongsiri \cite{srijuntongsiri_cci} proposed an iterative/subdivision hybrid algorithm for CCI that used a convergence test based on Kantorovich's theorem.  The algorithm, which was named \emph{the Kantorovich-Test Subdivision algorithm for Curve/Curve Intersection (KTS-CC)}, subdivides the parametric domains of the problems recursively while performing two tests on the subdomains until it is certain no additional subdivisions are required.  The first test uses a bounding volume of a subdomain to determine if the subdomain does not contain any solutions.  The second test is the convergence test based on Kantorovich's theorem, which tells us if Newton's method converges quadratically for the initial point in the current subdomain as well as whether it converges at all.  Once a subdomain passes the second test, KTS-CC invokes Newton's method to locate a solution, where quadratic convergence is ensured.  By appropriately combining the two tests and the subdivision scheme, the algorithm is guaranteed to find all of the intersection points between the two curves.  
The same idea was used for line/surface intersection \cite{srijuntongsiri_lsi} and for solving a system of polynomial equations with one degree of freedom \cite{srijuntongsiri_poly1free}.

The efficiency of KTS-CC depends on the total number of subdomains that need to be considered.  In particular, its convergence test requires choosing (or rather guessing) a test domain.  Choosing either too small or too large a test domain may cause the current subdomain to fail the test although Newton's method in fact converges quadratically to a solution, which results in unnecessary additional subdivisions.  KTS-CC always chooses the test domain to be 1.5 times the size of the current subdomain.  

In this article, we propose a modification to KTS-CC to adaptively adjust the size of the test domain according to what happens during the test of the parent subdomain.  In other words, if it appears that the parent test fails because the test domain is too large, we decrease the size for its children.  On the other hands, if it appears that the parent test fails because the domain is too small, we increase the size for its children.  We describe the adaptive algorithm in details in Section \ref{sec:algo}.  Experiment results comparing the efficiency of the adaptive algorithm to the original KTS-CC are in Section \ref{sec:experiment}.

Adaptive algorithms were proposed before for various problems such as adaptive conjugate gradients \cite{andrei,spillane,carson}, adaptive quadrature \cite{mckeeman,rice,Berntsen}, and adaptive simulated annealing \cite{ingber}.

\section{Curve/curve intersection problem}

Curve/curve intersection (CCI) is the problem of finding all of the intersection points between two space curves.  Among many conventional representations of space curves, we choose to use B\'{e}zier curves here.  Specifically, let $Z_{i,m}(t)$ denote the \emph{Bernstein polynomials}
\[
Z_{i,m}(t) = \frac{m!}{i!(m-i)!}(1-t)^{m-i}t^i.
\]
The two curves are represented by
\[
c^{(1)}(u) = \sum_{i=0}^m a_i Z_{i,m}(u), \mbox{\qquad} 0 \leq u \leq 1,
\]
where $a_i \in \mathbb{R}^3$ denote the coefficients, also known as the \emph{control points}, and
\[
c^{(2)}(v) = \sum_{j=0}^n a'_j Z_{j,n}(v), \mbox{\qquad} 0 \leq v \leq 1,
\]
where $a'_j \in \mathbb{R}^3$ are the control points.  The intersections between the two curves $c^{(1)}$ and $c^{(2)}$ are the solutions of
\begin{equation}
\label{eq1}
c^{(1)}(u)-c^{(2)}(v) = 0, \mbox{\qquad} 0 \leq u,v \leq 1,
\end{equation}
which is equivalent to
\begin{equation}
\label{maineq}
f(x) \equiv f(u,v) \equiv \sum_{i=0}^m \sum_{j=0}^n b_{ij} Z_{i,m}(u) Z_{j,n}(v) = 0, 
\end{equation}
where $0 \leq u,v \leq 1$ and $b_{ij} = a_i - a'_j \in \mathbb{R}^3$.  We denote $x = (u,v)^T$ and denote $f(x) \equiv f(u,v)$ for better readability below.  The system (\ref{maineq}) of three polynomials in two unknowns is the one our algorithm operates on.

\section{The theorem of Kantorovich}
\label{sec:kanto}
Denote the closed ball centered at $x \in \mathbb{R}^n$ with radius $r \in \mathbb{R}$, $r > 0$, by
\[
\bar{B}(x,r) = \{ y \in \mathbb{R}^n : \norm{y-x} \leq r \},
\]
and denote $B(x,r)$ as the interior of $\bar{B}(x,r)$.  Note that $x$ here refers to an arbitrary point in $\mathbb{R}^n$ and is different from $x$ in the previous section.  Kantorovich's theorem in affine invariant
form \cite{deuflhard,kantorovich} is

\begin{theorem}[Kantorovich, affine invariant form \cite{deuflhard,kantorovich}]
	\label{standardkantorovich}
	Let $f : D \subseteq \mathbb{R}^n \rightarrow \mathbb{R}^n$ be differentiable in
	the open convex set $D$. Assume that for some point $x^0 \in D$, the Jacobian $f'(x^0)$ is invertible.  Let $\eta$ be an upper bound
	\[
	\norm{f'(x^0)^{-1}f(x^0)} \leq \eta.
	\]
	Let there be a Lipschitz constant $\omega > 0$ for $f'(x^0)^{-1} f'$ such that
	\[
	\norm{f'(x^0)^{-1}(f'(x)-f'(y))} \leq \omega \cdot \norm{x-y} \textrm{ for all } x,y \in D.
	\]
	If $h = \eta\omega \leq 1/2$ and $\bar{B}(x^0,\rho_-) \subseteq D$, where
	\[
	\rho_- = \frac{1-\sqrt{1-2h}}{\omega},
	\]
	then $f$ has a zero $x^*$ in $\bar{B}(x^0,\rho_-)$. Moreover, this zero is the unique zero
	of $f$ in $(\bar{B}(x^0,\rho_-) \cup B(x^0,\rho_+)) \cap D$ where
	\[
	\rho_+ = \frac{1+\sqrt{1-2h}}{\omega}
	\]
	and the Newton iterates $x^k$ with
	\[
	x^{k+1} = x^k - f'(x^k)^{-1}f(x^k)
	\]
	are well-defined, remain in $\bar{B}(x^0,\rho_-)$, and converge to $x^*$. In addition,
	\begin{equation}
	\label{rapidconv}
	\norm{x^*-x^k} \leq \frac{\eta}{h}\left( \frac{(1-\sqrt{1-2h})^{2^k}}{2^k} \right), k = 0,1,2,\ldots
	\end{equation}
\end{theorem}
The point $x^0$ is said to be a \emph{fast starting point} if the sequence of Newton iterates starting from it
converges to a solution $x^*$ and (\ref{rapidconv}) is satisfied with $h \leq 1/4$,
which implies quadratic convergence of the iterates starting from $x^0$ \cite{srijuntongsiri_lsi}.
The Kantorovich's theorem also holds for complex functions \cite{farouki}.

\section{The adaptive Kantorovich-test subdivision algorithm for curve/curve intersection}
\label{sec:algo}

This section introduces the adaptive Kantorovich-test subdivision algorithm for curve/curve intersection algorithm (AKTS-CC).  
As we are interested in solutions of $f$ within the square $[0,1]^2$, and the closed ball $\bar B(x,r)$ defined in the infinity norm is a square, AKTS-CC uses the infinity norm for all of its norm computation.  It also maintains a list $S$ of \emph{explored regions} defined as parts of the domain $[0,1]^2$ guaranteed by Kantorovich's Theorem to contain only the zeros that have already been found. 

Let $f_i(x)$ denote the $i$th coordinate of the point $f(x)$ in three-dimensional space.  The  \emph{Kantorovich test} on a square $X = \bar B(x^0, r)$ is defined as the application of Kantorovich's theorem on the point $x^0$ to the functions $f_{ij}(x) = \left(f_i(x), f_j(x)\right)^T$ for each $\{i, j\} \in \left\{ \{1, 2\}, \{1, 3\}, \{2, 3\}\right\}$ using $\norm{f'_{ij}(x^0)^{-1}f_{ij}(x^0)}$ as $\eta$ in the statement of the theorem.  
For $\omega$, just as in the original KTS-CC, we use a more easily computed upper bound $\hat \omega \geq \omega$ as described in \cite{srijuntongsiri_cci} instead.  The square $X$ passes the Kantorovich test if there exists a pair $\{i,j\} \in \left\{ \{1, 2\}, \{1, 3\}, \{2, 3\}\right\}$ satisfying $\eta\hat\omega \leq 1/4$ and $\bar{B}(x^0, \rho_-) \subseteq D_{ij}$, where $D_{ij}$ is the chosen test domain for $f_{ij}$.  
We describe our adaptive scheme for determining $D_{ij}$ for each $X$ below as it involves the remaining parts of the algorithm.

If $X$ passes the Kantorovich test, $x^0$ is a fast starting point for $f_{ij}$ for the particular $\{i, j\}$ satisfying the conditions of the Kantorovich test, which means Newton's method starting from $x^0$ converges quadratically to a zero $x^*$ of $f_{ij}$.  We therefore perform Newton's method from $x^0$ to find $x^*$ and then evaluate $f(x^*)$ to see if $x^*$ is also a zero of $f$.  If it is, we have that 
\begin{equation}
\label{exp_comp}
X_E = \bar B(x^0, \rho_-) \cup \bar B(x^0, \rho_+),
\end{equation}
where $\rho_-$ and $\rho_+$ are as in the statement of Kantorovich's theorem, is an explored region associated with $x^*$.  Note that $X_E$ is defined in relation to $f_{ij}$ although it is also an explored region for $f$.

The other test used by our algorithm is the exclusion test.  For a given square $X$, let $\hat{f}_X$ be the Bernstein polynomial
that reparametrizes with $[0,1]^2$ the function defined by $f$ over
$X$.  
The square $X$ passes the \emph{exclusion test} if the convex hull of the control points of $\hat{f}_X$ excludes the origin.  As the exclusion test is the same as in KTS-CC, we refer the readers to \cite{srijuntongsiri_cci} for more details.

If a square $X = \bar B(x^0, r)$ fails the exclusion test, we then subdivide $X$ along both axes into four equal smaller squares $X_1$, $X_2$, $X_3$, and $X_4$ to be further investigated regardless of whether it passes the Kantorovich test (since there may be more than one zero in $X$ and passing the Kantorovich test guarantees converging to only one of them).  We choose the test domains for the subsquares $X_k$'s depending on how $X$ does on the Kantorovich test.  Specifically, suppose $D_{ij} = \bar B(x^0, \alpha_{ij} r )$.  
Let $D^k_{ij} = \bar B(x^0_k, \alpha^k_{ij} r')$ denote the test domain for $X_k$ with respect to $f_{ij}$, where $X_k = \bar B(x^0_k, r')$.
If $X$ passes the test (for any $f_{ij}$), we set $\alpha^k_{ij} = \alpha_{ij}$ for all $k$ and all pair $\{i,j\}$. 
If $X$ fails the test for $f_{ij}$ by having $\eta \hat \omega \leq 1/4$ but $\bar B(x^0, \rho_-) \not \subseteq D_{ij}$, we set $\alpha^k_{ij} = \alpha_{ij} + \epsilon$ for all $k$, where $\epsilon$ is a fixed constant.  The rationale is that it is possible that $x^0$ is in fact a fast starting point but we choose too small $D_{ij}$ for the test.  Specifically, although increasing $D_{ij}$ may increase $\hat \omega$, it may not cause $\eta \hat \omega$ to be greater than $1/4$ but make $\bar B(x^0, \rho_-)  \subseteq D_{ij}$.  On the other hand, if $X$ fails the test because $\eta \hat \omega > 1/4$, we set $\alpha^k_{ij} = \max \left(1, \alpha - \epsilon\right)$ for all $k$.  Similarly, $x^0$ may be a fast starting point but we choose too large $D_{ij}$ for the test (since decreasing the size of $D_{ij}$ may decrease $\hat \omega$ enough to satisfy $\eta \hat \omega \leq 1/4$).  We do not allow $\alpha^k_{ij}$ to be smaller than 1 as $D^k_{ij}$ should not be smaller than $X^k$ (otherwise, the Kantorovich test would not be able to detect a zero near a border of $X^k$).  For comparison, the original (nonadaptive) KTS-CC uses $\alpha_{ij} = 1.5$ for all $X$'s and all $\{i,j\}$'s.

Lastly, our algorithm maintains a first-in-first-out queue of areas in $[0,1]^2$ together with their three test domains that still need investigation. The details of AKTS-CC are as follows.
\begin{flushleft}
	\textbf{Algorithm AKTS-CC}:
\end{flushleft}
\begin{itemize}
	\item Initially, let the queue $Q$ contain only the square $[0,1]^2$ with its three test domains $D_{12} = D_{13} = D_{23} = \lbrack -0.25, 1.25\rbrack^2$.  Set $S = \emptyset$.
	\item While $Q$ is nonempty,
	\begin{enumerate}
		\item Let $X$ be the first square at the front of the queue $Q$.  Remove $X$ from $Q$.
		\item If $X \not\subseteq X_S$ for all $X_S \in S$,
		\begin{itemize}
			\item Perform the exclusion test on $X=\bar{B}(x^0,r)$.
			\item If $X$ fails the exclusion test,
			\begin{enumerate}
				\item Perform the Kantorovich test on $X$.
				\item If $X$ passes the Kantorovich test,
				\begin{enumerate}
					\item Perform Newton's method starting from $x^0$ to find a zero $x^*$ of $f_{ij}$, where $\{i, j\}$ is the pair with which $X$ passes the Kantorovich test.
					\item If $f(x^*) = 0$ and $x^* \not\in X_S$ for any $X_S \in S$ (i.e., $x^*$ is a zero of $f$ and has not been found previously),
					\begin{itemize}
						\item Record $x^*$ as one of the intersections.
						\item Compute the new associated explored region $X_E$ according to (\ref{exp_comp}).
						\item Set $S = S \cup X_E$.
					\end{itemize}
				\end{enumerate}
				\item Subdivide $X$ along both axes into $4$ equal smaller squares. Add these squares to the end of $Q$ together with their test domains defined as described above.
			\end{enumerate}
		\end{itemize}
	\end{enumerate}
\end{itemize}

\section{Computational results}
\label{sec:experiment}

\begin{figure}
	\centering
	\includegraphics[width=1.0\textwidth]{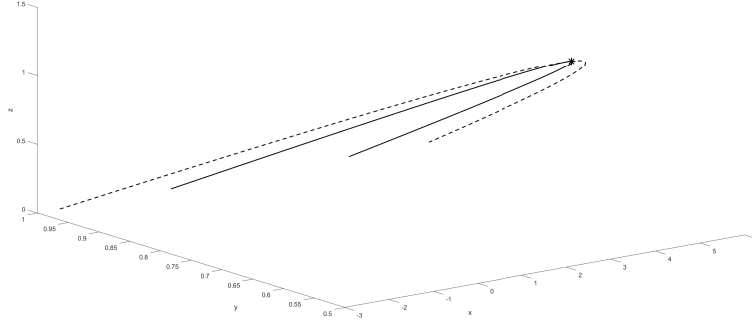}
	\caption{Curves of test problem 1 and their intersection.  The intersection is shown with *.}
	\label{fig:p1}
\end{figure}

\begin{figure}
	\centering
	\includegraphics[width=1.0\textwidth]{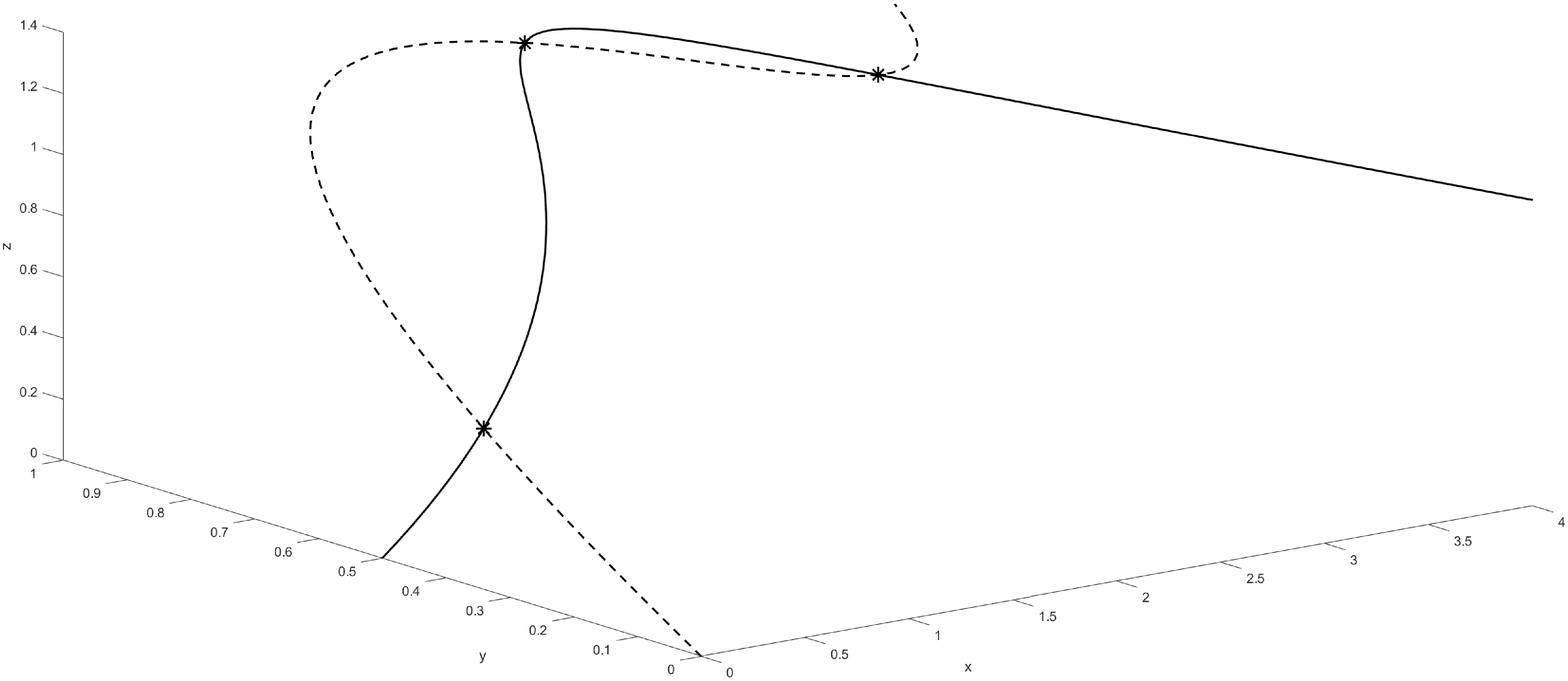}
	\caption{Curves of test problem 2 and their intersections.  The intersections are shown with *.}
	\label{fig:p2}
\end{figure}

\begin{figure}
	\centering
	\includegraphics[width=1.0\textwidth]{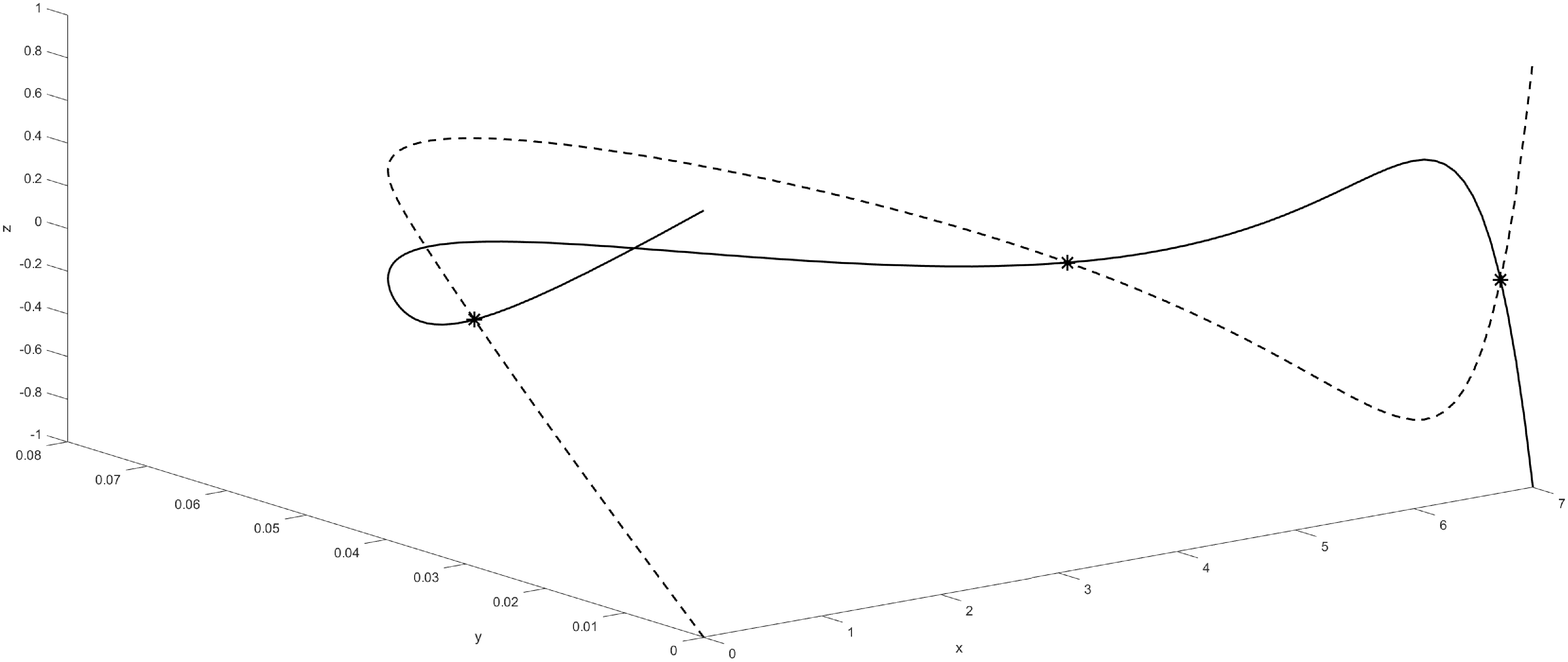}
	\caption{Curves of test problem 3 and their intersections.  The intersections are shown with *.}
	\label{fig:p3}
\end{figure}

\begin{figure}
	\centering
	\includegraphics[width=1.0\textwidth]{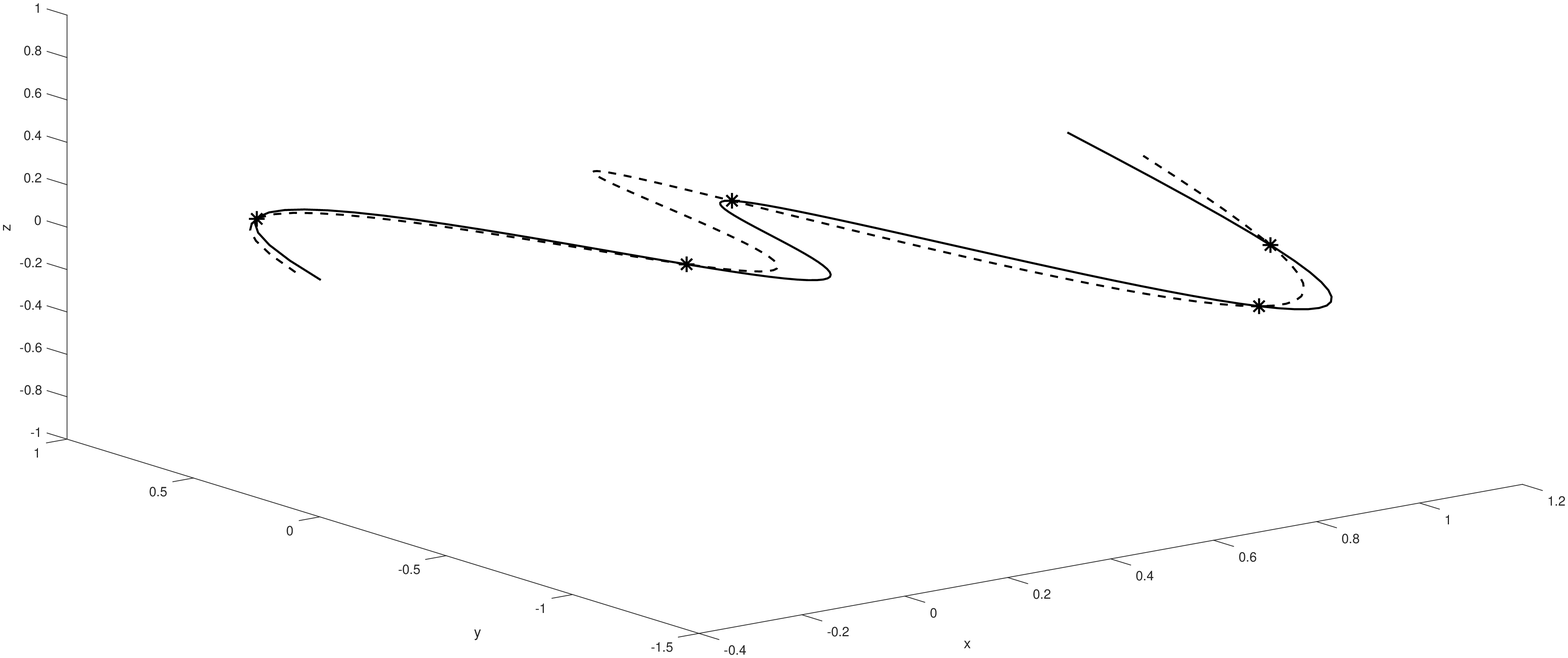}
	\caption{Curves of test problem 4 and their intersections.  The intersections are shown with *.}
	\label{fig:p4}
\end{figure}

\begin{figure}
	\centering
	\includegraphics[width=1.0\textwidth]{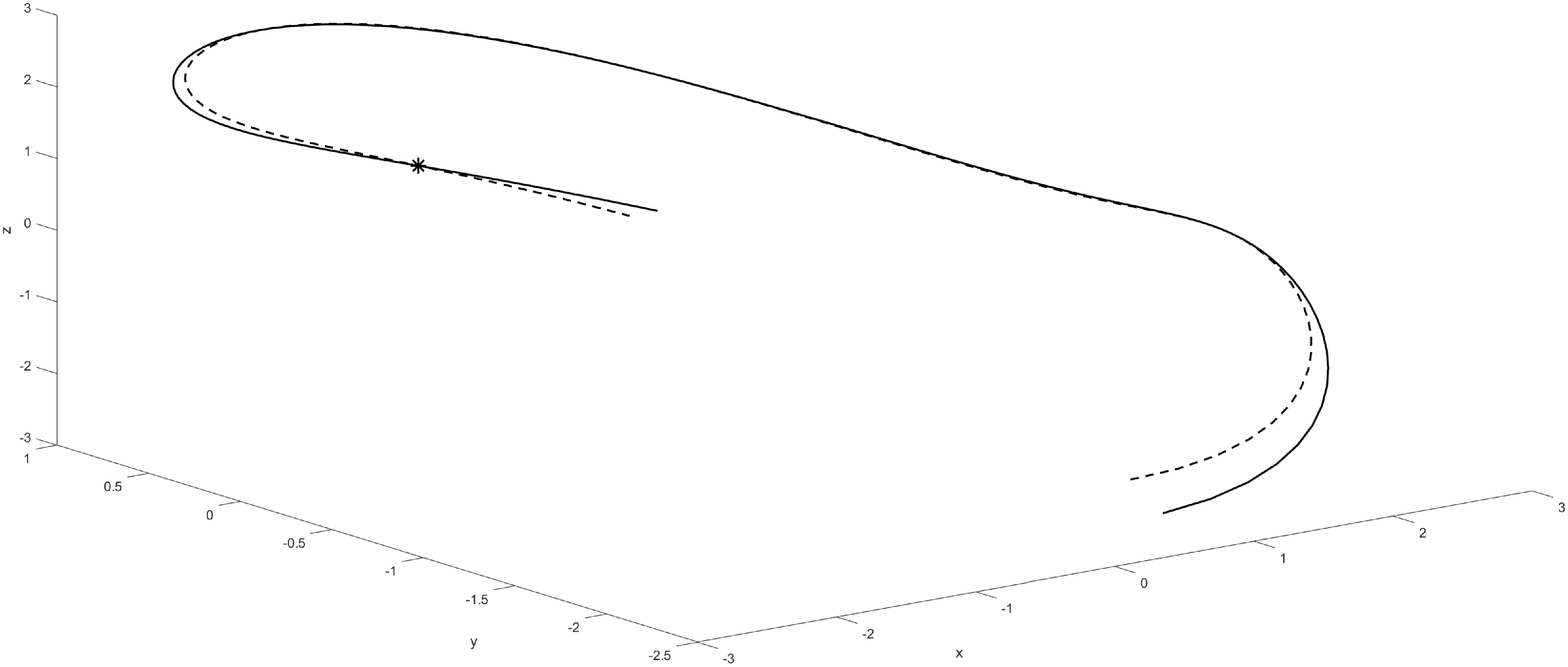}
	\caption{Curves of test problem 5 and their intersection.  The intersection is shown with *.}
	\label{fig:p5}
\end{figure}

\begin{figure}
	\centering
	\includegraphics[width=1.0\textwidth]{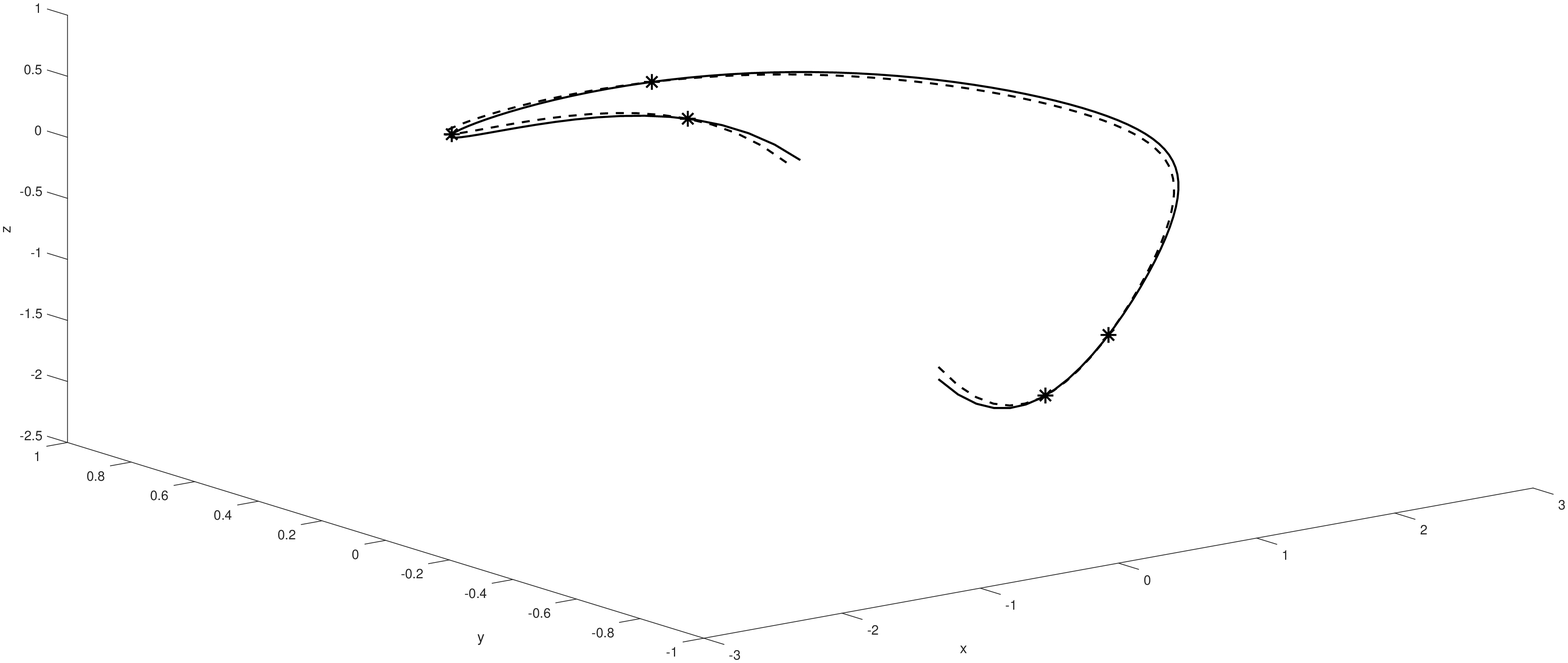}
	\caption{Curves of test problem 6 and their intersections.  The intersections are shown with *.}
	\label{fig:p6}
\end{figure}

\begin{figure}
	\centering
	\includegraphics[width=1.0\textwidth]{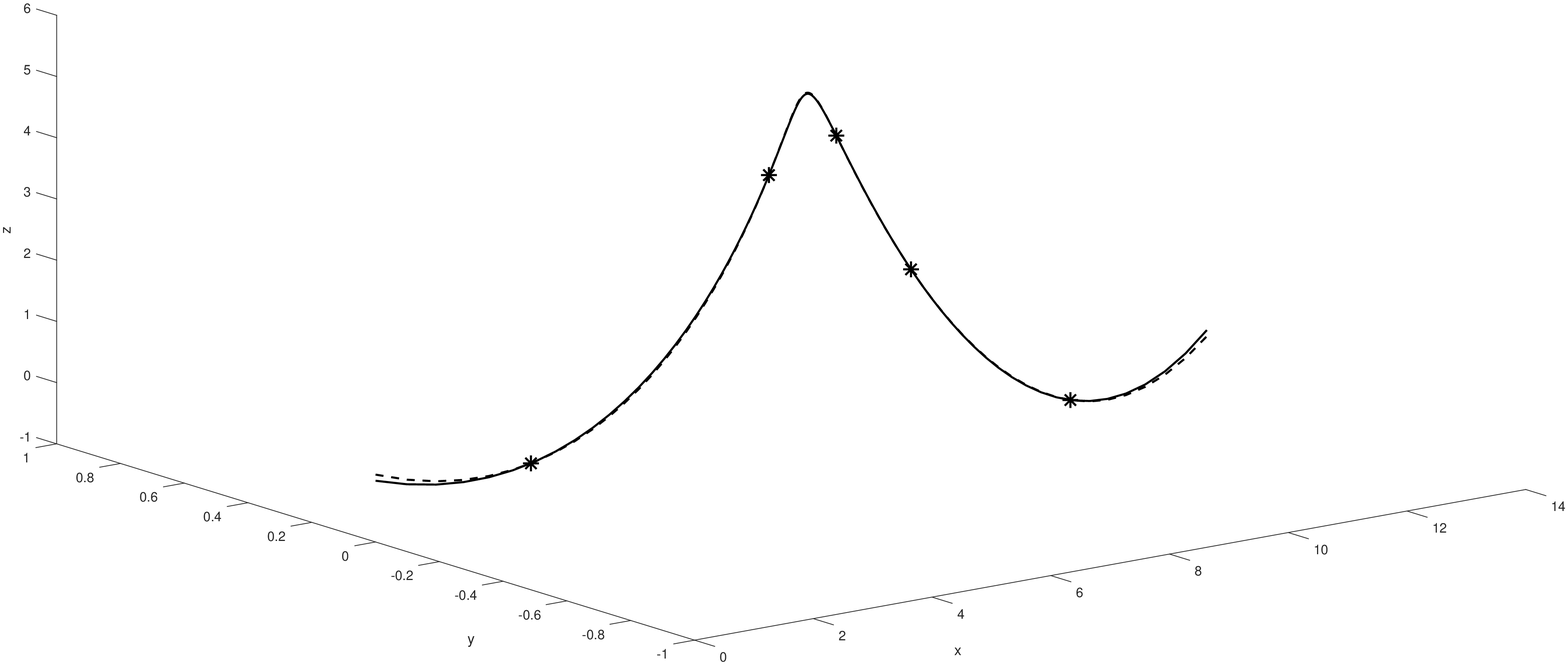}
	\caption{Curves of test problem 7 and their intersections.  The intersections are shown with *.}
	\label{fig:p7}
\end{figure}

\begin{figure}
	\centering
	\includegraphics[width=1.0\textwidth]{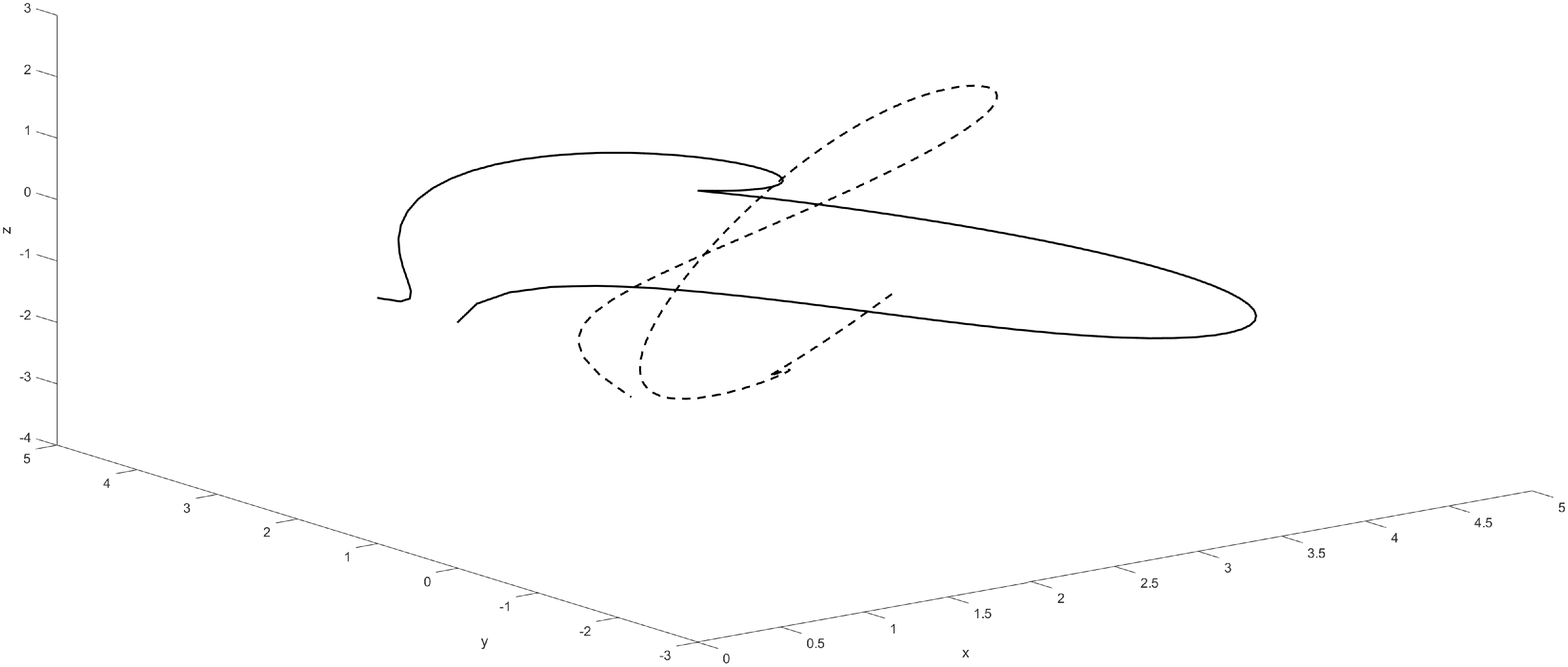}
	\caption{Curves of test problem 8, which do not intersect.}
	\label{fig:p8}
\end{figure}

We implemented the proposed adaptive algorithm AKTS-CC in Matlab and compared its efficiency with varying values of $\epsilon$ against KTS-CC on a number of test problems.  Our implementation uses the reparametrization algorithm presented in \cite{srijuntongsiri_basis}.  The experiments were performed using tolerance of $10^{-7}$ for Newton's method parts of both algorithms.  The test problems and their intersections computed by AKTS-CC are shown in Figures \ref{fig:p1} to \ref{fig:p8}.  
Table \ref{table:result} compares the efficiency of AKTS-CC and KTS-CC on the eight test problems.  It reports the degrees of the two curves and the total number of squares examined by AKTS-CC for different values of $\epsilon$ and KTS-CC during the entire computations.  Since the number of operations per While-loop iteration in AKTS-CC is only a small constant larger than in KTS-CC, the efficiency of the two algorithms depend on the total number of iterations, which is the same as the number of squares examined.

\begin{table}
	\centering
	\caption{The total number of squares examined by AKTS-CC for varying values of $\epsilon$ and KTS-CC on the test problems shown in Figures \ref{fig:p1}--\ref{fig:p8}.}
	\label{table:result}
	\begin{tabular}{rrrrrrrr}
		\hline %
		\multicolumn{1}{c}{\multirow{2}{*}{Fig.}} & Degrees & \multicolumn{5}{c}{AKTS-CC}  & \multirow{2}{*}{KTS-CC} \\
		\cline{3-7} 
		& $\left((m, n)\right)$ & \multicolumn{1}{c}{$\epsilon = 0.01$} & \multicolumn{1}{c}{$\epsilon = 0.05$} & \multicolumn{1}{c}{$\epsilon = 0.1$} & \multicolumn{1}{c}{$\epsilon = 0.15$} & \multicolumn{1}{c}{$\epsilon = 0.2$} \\
\hline
		\ref{fig:p1} & $(2, 2)$ & 5&5&5&5&5&5  \\
		\ref{fig:p2} & $(3, 3)$ & 41&41&41&41&41&41  \\
		\ref{fig:p3} & $(7, 7)$ & 45 & 45 &45&45&45&49  \\
		\ref{fig:p4} & $(8, 8)$ & 149&145&145&145&145&149  \\
		\ref{fig:p5} & $(9, 9)$ & 85&81&81&81&81&85  \\
		\ref{fig:p6} & $(9, 9)$ & 253&253&253&253&253&253  \\
		\ref{fig:p7} & $(9, 9)$ & 313&313&313&313&313&313  \\
		\ref{fig:p8} & $(13, 12)$ & 29&29&29&29&29&29  \\
\hline
	\end{tabular}
\end{table}

The results show that AKTS-CC is slightly more efficient than KTS-CC in three of the eight test problems and is as efficient as KTS-CC in the five remaining ones.  Additionally, AKTS-CC saves just four squares in those cases.  On the other hand, as AKTS-CC differ from KTS-CC only on the choices of the test domains for the Kantorovich test, which matters when the current square contains a zero, we do not expect improvement on squares not containing any zeros in any case.  

\section{Conclusion}

We propose Algorithm AKTS-CC that is a modification of KTS-CC to adaptively determine the test domains for the Kantorovich test based on the result of the same test on the parent square.  The test domains for $f_{12}$, $f_{13}$, and $f_{23}$ for the same square are determined independently.  The algorithm was implemented in Matlab and was shown to be marginally more efficient than KTS-CC for some test cases and equally efficient for others.  Future investigations on different adaptive schemes may result in larger improvement of the efficiency of the algorithm.

\section*{Acknowledgement}
The author gratefully acknowledges the financial support provided by Thammasat University Research Fund under the TU Research Scholar, Contract No. TP 2/24/2560.

\bibliography{allbib}

\end{document}